\documentclass{ifacconf}
\usepackage[centertags]{amsmath}
\usepackage{graphics}
\usepackage{amscd}
\usepackage{amsfonts}
\usepackage{amssymb}
\usepackage{helvet}
\usepackage{rotating}
\usepackage{epsfig}
\usepackage[latin1]{inputenc}
\usepackage{placeins}
\usepackage{pst-all}
\usepackage{relsize}
\usepackage{longtable}
\usepackage{natbib}

\usepackage{ifpdf}
\ifpdf%
 \usepackage{pdftricks}
 \begin{psinputs}
    \usepackage{pstricks}
    \usepackage{pstricks-add}
    \usepackage{pst-node}
    \usepackage{pst-grad}
    \usepackage{pst-text}
  \end{psinputs}
\else
  \usepackage{pstricks}
  \usepackage{pstricks-add}
  \usepackage{pst-node}
  \usepackage{pst-grad}
  \usepackage{pst-text}
  \newenvironment{pdfpic}{}{}
\fi

\newtheorem{theorem}{Theorem}[section]

\newtheorem{proposition}[theorem]{Proposition}
\newtheorem{definition}[theorem]{Definition}
\newtheorem{assumption}[theorem]{Assumption}

\newcommand{\R}{\mathbb{R}}
\newcommand{\N}{\mathbb{N}}

\newcommand{\X}{\mathbb{X}}
\newcommand{\U}{\mathbb{U}}
\newcommand{\XX}{X}
\newcommand{\UU}{U}
\newcommand{\cI}{\mathcal{I}}

\def\argmin{\mathop{\rm argmin}}

\usepackage{framed}
\newenvironment{fshaded}{%
\MakeFramed {\FrameRestore}}%
{\endMakeFramed}
\definecolor{shadecolor}{rgb}{1.,1.,1.}%
\definecolor{framecolor}{rgb}{.0,.0,.0}%
\begin{document}
%%%%%%%%%%%%%%%%%%%%%%%%%%%%%%%%%%%%%%%%%%%%%%%%%%%%%%%%%%%%%%%%%%%%%%%%%%%%%%%%
\begin{frontmatter}

\title{Stability of Constrained Adaptive Model Predictive Control Algorithms\thanksref{footnoteinfo}}

\thanks[footnoteinfo]{This work was supported by DFG Grant Gr1569/12 within the Priority Research Program 1305 and the Leopoldina Fellowship Programme LPDS 2009-36.}

\author[First]{T. Jahn} 
\author[Second]{J. Pannek} 

\address[First]{University of Bayreuth, 95440 Bayreuth, Germany\newline (e-mail: thomas.jahn@uni-bayreuth.de)}
\address[Second]{Curtin University of Technology, Perth, 6845 WA, Australia\newline (e-mail: juergen.pannek@googlemail.com)}

\begin{abstract}
Recently, suboptimality estimates for model predictive controllers (MPC) have been derived for the case without additional stabilizing endpoint constraints or a Lyapunov function type endpoint weight. The proposed methods yield \textit{a posteriori} and \textit{a priori} estimates of the degree of suboptimality with respect to the infinite horizon optimal control and can be evaluated at runtime of the MPC algorithm. Our aim is to design automatic adaptation strategies of the optimization horizon in order to guarantee stability and a predefined degree of suboptimality for the closed loop solution. Here, we present a stability proof for an arbitrary adaptation scheme and state a simple shortening and prolongation strategy which can be used for adapting the optimization horizon.
\end{abstract}


\end{frontmatter}
%%%%%%%%%%%%%%%%%%%%%%%%%%%%%%%%%%%%%%%%%%%%%%%%%%%%%%%%%%%%%%%%%%%%%%%%%%%%%%%%

%%%%%%%%%%%%%%%%%%%%%%%%%%%%%%%%%%%%%%%%%%%%%%%%%%%%%%%%%%%%%%%%%%%%%%%%%%%%%%%%
\section{Introduction}
\label{Section:Introduction}
%%%%%%%%%%%%%%%%%%%%%%%%%%%%%%%%%%%%%%%%%%%%%%%%%%%%%%%%%%%%%%%%%%%%%%%%%%%%%%%%

Nowadays, model predictive controllers (MPC), sometimes also called receding horizon controllers (RHC), are used in a variety of industrial applications, cf. \citet{QB2003}. As shown in \citet{AZ2000}, \citet{MRRS2000} and \citet{RM2009}, theory for such controllers is also widely understood both for linear and nonlinear systems. The control method itself deals with the problem of approximately solving an infinite horizon optimal control problem which is computationally intractable in general. Reasons for its success are on the one hand its capability to directly incorporate constraints depending on the states and inputs of the underlying process. On the other hand, the fundamental steps of this method are very simple: First, a solution of a finite horizon optimal control problem is computed for a given initial value. In a second step, the first element of the resulting control is implemented at the plant and in the last step, the finite horizon is shifted forward in time. As a consequence, the method is iteratively applicable and reveals the control to be a static state feedback.

Unfortunately, stability of solution of the infinite horizon problem may be lost due to considering only finite horizons. Over the last two decades, several solutions have been proposed to cope with this issue, see, e.g., \citet{KG1988}, \citet{CA1998} and \citet{GR2008}. All these approaches require the horizon to be sufficiently long and computing the minimal required horizon length is computationally demanding. However, the horizon needs to be chosen as a worst case scenario which is usually needed to cope with small regions of state space only. Our aim in this work is to develop online applicable adaptation strategies for the horizon length which guarantee stability of the closed loop. In particular, we follow the approach of \citet{GP2009} where different suboptimality estimates have been developed to measure the performance of the model predictive controller. Based on these estimates, we propose a simple technique to locally fit the horizon to the control task, the current state of the system and also to the MPC internal information. Due to the change of the structure of the controller, however, known stability proofs and suboptimality results cannot be applied. To cover these issues, we present a stability result for MPC with varying optimization horizons using mild additional conditions. To some extend adaptation strategies of the horizon are known in the literature, see e.g. \citet{FLJ2006} and \citet{VBRSB2008}, which are heuristics based on insight of the specific problem but have shown to be applicable in an adaptive model predictive control setting. In contrast to that, our approach can be proven rigorously and doesnot require any insight into the process under consideration (note that different to our intention the term adaptive model predictive control is also used to incorporate model uncertainties, see, e.g., \citet{MM1993} and \citet{AG2009}).

The paper is organized as follows: In Section \ref{Section:Setting} we describe and motivate the problem setup. Section \ref{Section:Stability for Standard NMPC} deals with the a posteriori and a priori suboptimalty estimates which will be the foundation of our analysis. In the following Section \ref{Section:Stability under Adaptation}, we show how the stated stability results and estimates can be extended to the case of varying optimization horizons. Thereafter, we state a simple shortening and prolongation strategy based on the suboptimality estimates given in Section \ref{Section:A simple Adaptation Strategy}. In order to show the applicability and effectiveness of our approach, Section \ref{Section:Numerical Results} contains a numerical example of the adaptive MPC approach. The final Section \ref{Section:Conclusion} concludes the paper and points out directions of future research.

%%%%%%%%%%%%%%%%%%%%%%%%%%%%%%%%%%%%%%%%%%%%%%%%%%%%%%%%%%%%%%%%%%%%%%%%%%%%%%%%
\section{Setting}
\label{Section:Setting}
%%%%%%%%%%%%%%%%%%%%%%%%%%%%%%%%%%%%%%%%%%%%%%%%%%%%%%%%%%%%%%%%%%%%%%%%%%%%%%%%

Within this work we analyze nonlinear discrete time control systems of the form
\begin{align}
	\label{Setup:nonlinear discrete time system}
	x(i + 1) = f(x(i), u(i)), \quad x(0) = x_0
\end{align}
with $x(i) \in \X \subset \XX$ and $u(i) \in \U \subset \UU$ for $i \in \N_0$. For the considered systems the state space $\XX$ and the control value space $\UU$ are arbitrary metric spaces. Hence, all following results also apply to the discrete time dynamics induced by a sampled infinite dimensional system, cf.\ \citet{IK2002} or \citet{AGW2010a}. Here, we denote the space of control sequences $u: \N_0 \rightarrow \U$ by $\U^{\N_0}$ and the solution trajectory for given control $u \in \U^{\N_0}$ by $x_u(\cdot)$. Additionally, the sets $\X$ and $\U$ incorporate possible restrictions on the state and control respectively.

In the following, we aim at finding a static state feedback $u = \mu(x) \in \U^{\N_0}$ for a given control system \eqref{Setup:nonlinear discrete time system} which minimizes the infinite horizon cost functional
\begin{align}
	\label{Setup:infinite cost functional}
	J_\infty (x_0, u) = \sum\limits_{i = 0}^\infty l(x_u(i), u(i))
\end{align}
with stage cost $l: \XX \times \UU \rightarrow \R_0^+$. The corresponding optimal value function is denoted by $V_\infty(x_0) = \inf_{u \in \U^{\N_0}} J_\infty(x_0, u)$ and throughout this paper we assume that the minimum with respect to $u \in \U^{\N_0}$ is attained. The optimal value function $V_\infty(\cdot)$ can be used to define the infinite horizon feedback law
\begin{align}
	\label{Setup:infinite control}
	\mu(x_u(i)) := \argmin_{u \in \U} \left\{ V_\infty(x_u(i + 1)) + l(x_u(i), u) \right\}
\end{align}
for which one can show optimality using Bellman's optimality principle. Since the computation of the desired control law requires the solution of a Hamilton--Jacobi--Bellman equation, we use a model predictive control approach in order to avoid the problem of solving an infinite horizon optimal control problem. The fundamental idea of such a model predictive controller is simple and consists of three steps which are repeated at every discrete time instant during the process run: First, an optimal control for the problem on a finite horizon $[0, N]$ is computed given the most recent known state of the system $x_0$. Then, the first control element is implemented at the plant and in the third step the entire optimal control problem considered in the first step is shifted forward in time by one discrete time instant which allows for iteratively repeating this process. In the literature this method is also termed receding horizon, see, e.g., \citet{MRRS2000}.

In contrast to the infinite horizon optimal control \eqref{Setup:infinite control}, the problem in the second step is to minimize the truncated cost functional on a finite horizon
\begin{align}
	\label{Setup:finite cost functional}
	J_N(x_0, u) = \sum\limits_{k = 0}^{N - 1} l(x_u(k, x_0), u(k)).
\end{align}
The truncated horizon defines the set of discrete time instances $\cI := \{ 0, \ldots, N - 1 \}$. Here, we assume the first instant to be denoted by zero for each optimal control problem within the MPC problem. In particular, we focus on the implementation of a constrained model predictive controller without additional stabilizing endpoint constraints or a Lyapunov function type endpoint weight, see, e.g., \citet{KG1988} and \citet{CA1998}, respectively.

Throughout this work, we denote the closed loop solution at time instant $i$ by $x(i)$ while $x_u(\cdot, x_0)$ denotes the open loop trajectory of the prediction. Moreover, we use the abbreviations
\begin{align}
	\label{Setup:open loop control}
	u_N(\cdot, x_0) & = \argmin_{u \in \U^N} J_N(x_0, u) \\
	u_N(x_0) & = u_N(0, x_0) \nonumber
\end{align}
for the minimizing open loop control sequence of the reduced cost functional and its first element respectively. We call $V_N(x_0) = \min_{u \in \U^{N}} J_N(x_0, u)$ the optimal value function of the finite cost functional \eqref{Setup:finite cost functional} and, for notational purposes, we use $u_N(i, x_0)$ to represent the $i$-th control value within the open loop control sequence corresponding to the initial value $x_0$ when it is necessary to distinguish between two or more different open loop controls. Hence, if the initial value $x_{u_N}(0, x_0) = x_0$ is given, then the open loop control \eqref{Setup:open loop control} induces the open loop solution
\begin{align}
	\label{Setup:open loop solution}
	x_{u_N}(k + 1, x_0) = f\left( x_{u_N}(k, x_0), u_N(k, x_0) \right)
\end{align}
for all time instances $k$ on the optimization horizon $\cI \setminus \{N\}$. Similarly to \eqref{Setup:infinite control}, the closed loop control can be defined as
\begin{align}
	\label{Setup:closed loop control}
	\mu_N(x(i)) := \argmin_{u \in \U} \left\{ V_{N - 1}(x(i + 1)) + l(x(i), u) \right\}
\end{align}
and the corresponding closed loop system is given by
\begin{align}
	\label{Setup:closed loop solution}
	x(i + 1) = f\left( x(i), \mu_N(x(i)) \right)
\end{align}
for all $i \in \N_0$.

Note that due to the truncation of the infinite horizon cost functional \eqref{Setup:infinite cost functional} to the finite MPC cost functional \eqref{Setup:finite cost functional}, stability and optimality properties of the closed loop solution \eqref{Setup:closed loop control}, \eqref{Setup:closed loop solution} induced by the infinite horizon optimal control \eqref{Setup:infinite control} are not preserved in general.

Here, our aim is to show that in order to guarantee stability of the closed loop \eqref{Setup:closed loop control}, \eqref{Setup:closed loop solution} for any initial value $x \in \X$, the requirement of considering the worst case optimization horizon $N$ for all initial values $x \in \X$ can be weakened. Additionally, the resulting closed loop trajectory satisfies locally a predefined degree of suboptimality compared to the infinite horizon solution \eqref{Setup:nonlinear discrete time system}, \eqref{Setup:infinite control}.

%%%%%%%%%%%%%%%%%%%%%%%%%%%%%%%%%%%%%%%%%%%%%%%%%%%%%%%%%%%%%%%%%%%%%%%%%%%%%%%%
\section{Stability for Standard NMPC}
\label{Section:Stability for Standard NMPC}
%%%%%%%%%%%%%%%%%%%%%%%%%%%%%%%%%%%%%%%%%%%%%%%%%%%%%%%%%%%%%%%%%%%%%%%%%%%%%%%%

The measure of suboptimality we consider in the following is the difference between the infinite horizon cost induced by the MPC law $\mu_N(\cdot)$, that is
\begin{align}
	\label{Suboptimality:eq:value function mu_N}
	V_\infty^{\mu_N} (x_0) := \sum\limits_{i = 0}^\infty l\left( x(i), \mu_N(x(i)) \right),
\end{align}
and the finite horizon cost $V_N(\cdot)$ or the infinite horizon optimal value function $V_\infty(\cdot)$. In particular, the latter give us estimates on the degree of suboptimality of the controller $\mu_N(\cdot)$ of the MPC process. For this purpose, we make extensive use of the suboptimality estimates derived in \citet{GP2009}.

\begin{proposition}[A posteriori Estimate]\label{Suboptimality:prop:trajectory a posteriori estimate}
	Consider a feedback law $\mu_N: \X \rightarrow \U$ and its associated trajectory $x(\cdot)$ according to \eqref{Setup:closed loop solution} with initial value $x(0) = x_0 \in \X$. If there exists a function $V_N: \X \rightarrow \R_0^+$ satisfying
	\begin{align}
		\label{Suboptimality:prop:trajectory a posteriori estimate:eq1}
		V_N(x(i)) \geq V_N(x(i + 1)) + \alpha l(x(i), \mu_N(x(i)))
	\end{align}
	for some $\alpha \in (0, 1]$ and all $i \in \N_0$, then
	\begin{align}
		\label{Suboptimality:prop:trajectory a posteriori estimate:eq2}
		\alpha V_{\infty}(x(i)) \leq \alpha V_{\infty}^{\mu_N}(x(i)) \leq V_N(x(i)) \leq V_\infty(x(i))
	\end{align}
	holds for all $i \in \N_0$.
\end{proposition}

Since all values in \eqref{Suboptimality:prop:trajectory a posteriori estimate:eq1} are computed throughout the NMPC process, $\alpha$ can be easily computed online along the closed loop trajectory. Thus, \eqref{Suboptimality:prop:trajectory a posteriori estimate:eq1} yields a computationally feasible and numerically cheap way to estimate the degree of suboptimality of the trajectory.

Due to the fact that $V_N(x(i + 1))$ in \eqref{Suboptimality:prop:trajectory a posteriori estimate:eq1} is unknown at runtime, Proposition \eqref{Suboptimality:prop:trajectory a posteriori estimate} yields an a posteriori estimator. However, we can also utilize a more conservative a priori estimate if we assume the following:

\begin{assumption}\label{Suboptimality:ass:apriori}
	For given $N$, $N_0 \in \N$, $N \geq N_0 \geq 2$, there exists a constant $\gamma > 0$ such that for the open loop solution $x_{u_N}(i, x(i))$ given by \eqref{Setup:open loop solution} the inequalities
	\begin{align*}
		& \frac{V_{N_0}(x_{u_N}(N - N_0, x(i)))}{\gamma + 1} \leq \\
		& \leq \max_{j = 2, \ldots, N_0} l(x_{u_N}(N - j, x(i)), \mu_{j - 1}(x_{u_N}(N - j, x(i)))) \\
		& \frac{V_k(x_{u_N}(N - k, x(i)))}{\gamma + 1} \leq \\
		& \leq l(x_{u_N}(N - k, x(i)), \mu_k(x_{u_N}(N - k, x(i))))
	\end{align*}
	hold for all $k \in \{N_0 + 1, \ldots, N\}$ and all $i \in \N_0$.
\end{assumption}

\begin{theorem}[A priori Estimate]\label{Suboptimality:thm:apriori}
	Consider $\gamma > 0$ and $N$, $N_0 \in \N$, $N \geq N_0$ such that $(\gamma + 1)^{N - N_0} > \gamma^{N - N_0 + 2}$ holds. If Assumption \ref{Suboptimality:ass:apriori} is fulfilled for these $\gamma$, $N$ and $N_0$, then the estimate \eqref{Suboptimality:prop:trajectory a posteriori estimate:eq2} holds for all $i \in \N_0$ where
	\begin{align}
		\label{Suboptimality:thm:apriori:eq1}
		\alpha := \frac{(\gamma + 1)^{N - N_0} - \gamma^{N - N_0 + 2}}{(\gamma + 1)^{N - N_0}}.
	\end{align}
\end{theorem}

Note that we cannot expect the relaxed Lyapunov inequality \eqref{Suboptimality:prop:trajectory a posteriori estimate:eq1} or Assumption \ref{Suboptimality:ass:apriori} to hold in practice. In many cases the discrete time system \eqref{Setup:nonlinear discrete time system} is obtained from a discretization of a continuous time system, e.g. sampling with zero order hold, see \citet{NT2004}. Hence, even if the continuous time system is stabilizable to a setpoint $x^*$ and no numerical errors occur during optimization and integration, the corresponding sampled--data system is most likely practically stabilizable at $x^*$ only. However, suboptimality results can be extended to cover the case of practical stability as well, see \citet{GR2008} and \citet{GP2009}. Since extending the stability results we will present now to cover the practical case can be done analogously, see \citet{P2010}, we restrict ourselves to the case of asymptotic stability for simplicity of exposition.

%%%%%%%%%%%%%%%%%%%%%%%%%%%%%%%%%%%%%%%%%%%%%%%%%%%%%%%%%%%%%%%%%%%%%%%%%%%%%%%%
\section{Stability under Adaptation}
\label{Section:Stability under Adaptation}
%%%%%%%%%%%%%%%%%%%%%%%%%%%%%%%%%%%%%%%%%%%%%%%%%%%%%%%%%%%%%%%%%%%%%%%%%%%%%%%%

As stated at the end of Section \ref{Section:Setting}, we aim at weakening the worst case nature of the optimization horizon $N$. Here, one has to keep in mind that if a model predictive controller shall be designed for a given application, then stability of the resulting closed loop \eqref{Setup:closed loop solution} needs to be guaranteed for the entire working range $\X$. In practice, this may lead to very large optimization horizons $N$. Yet, most points visited by the closed loop \eqref{Setup:closed loop solution} we do not require such a large optimization horizon in order to guarantee stability.

Here, we focus on locally guaranteeing a decrease of the cost function for each step of the MPC process and modify the horizon length $N$ to fulfill this task. Similar to the suboptimality results from Section \ref{Section:Stability for Standard NMPC}, we want to measure this decrease in terms of the running cost $l(\cdot, \cdot)$ such that a given suboptimality bound $\overline{\alpha} \in (0, 1)$ is locally satisfied.

Since we are now dealing with varying optimization horizons, we intuitively extend our notation from Section \ref{Section:Stability for Standard NMPC} by adding the used optimization horizon as an argument, i.e. $\alpha(N)$ denotes the suboptimality degree $\alpha$ with horizon $N$. Moreover, since the resulting closed loop control now depends on a sequence $(N_i)_{i \in \N}$ we denote such a control law by $\mu_{(N_i)}$.

An abstract adaptive MPC algorithm which locally accomplishes the task of guaranteeing a decrease in the cost function is the following:

\begin{fshaded}
	\begin{itemize}
		\item[(1)] Given $x(i)$ and $N_i$ do
		\begin{itemize}
			\item[(1a)] Compute optimal control on horizon $N_i$
			\item[(1b)] Compute suboptimality degree $\alpha(N_i)$
			\item[(1c)] If $\alpha(N_i) \geq \overline{\alpha}$: Call shortening strategy for $N_i$ \\
			Else: Call prolongation strategy for $N_i$
		\end{itemize}
		while $\alpha(N_i) \leq \overline{\alpha}$
		\item[(2)] Implement the first control component $\mu_{N_i}(x(i)) := u(0, x(i))$
		\item[(3)] Set $i := i + 1$ and shift the optimization horizon forward in time
	\end{itemize}
\end{fshaded}

In this context, we distinguish the following degrees of suboptimality:

\begin{definition}[Suboptimality Degree]\label{Stability under Adaptation:def:suboptimality degree}
	(i) Given a set $\X$, then we call $\alpha := \max \{ \alpha \mid \text{\eqref{Suboptimality:prop:trajectory a posteriori estimate:eq1} holds $\forall x(n) = x \in \X$} \}$ the {\it global suboptimality degree}. \\
	(ii) Given a point $x \in \X$, then  we call $\alpha := \max \{ \alpha \mid \text{\eqref{Suboptimality:prop:trajectory a posteriori estimate:eq1} holds for $x(n) = x$} \}$ the {\it local suboptimality degree}. \\
	(iii) Given a closed loop trajectory $x(\cdot)$ we call $\alpha := \max \{ \alpha \mid \text{\eqref{Suboptimality:prop:trajectory a posteriori estimate:eq1} holds $\forall n \in \N_0$} \}$ the {\it closed loop suboptimality degree}.
\end{definition}

The problem which we are facing for such an adaptive MPC algorithm is the fact that none of the existing stability proofs, see, e.g., \citet{KG1988}, \citet{CA1998}, \citet{GMTT2005}, \citet{JH2005}, \citet{GP2009} %
%, \citet{Graichen2010}
and \citet{GPSW2010}, can be applied in this context since these results assume $N$ to be constant while here the optimization horizon $N_i$ may change in every step of the MPC algorithm.

The major obstacle to apply the idea of Proposition \ref{Suboptimality:prop:trajectory a posteriori estimate} in the context of varying optimization horizons $N$ is the lack of a common Lyapunov function along the closed loop. To compensate for this deficiency, we make the following mild assumption:

\begin{assumption}\label{Stability under Adaptation:ass:enhanced stabilizing}
	Given an initial value $x \in \X$ and a horizon length $N < \infty$ such that $\mu_N(\cdot)$ guarantees local suboptimality degree $\alpha(N) \geq \overline{\alpha}$, $\overline{\alpha} \in (0, 1)$, we assume that for $\widetilde{N} \geq N$, $\widetilde{N} < \infty$, there exist constants $C_l, C_\alpha > 0$ such that the inequalities
	\begin{align*}
		C_l l(x, \mu_{N}(x)) & \leq l(x, \mu_{\widetilde{N}}(x)) \frac{V_{\widetilde{N}}(x) - V_{\widetilde{N}}(f(x, \mu_N(x)))}{V_{\widetilde{N}}(x) - V_{\widetilde{N}}(f(x, \mu_{\widetilde{N}}(x))} \\
		C_\alpha \alpha(N) & \leq \alpha(\widetilde{N})
	\end{align*}
	hold where $\alpha(\widetilde{N})$ is the local suboptimality degree of the controller $\mu_{\widetilde{N}}(\cdot)$ corresponding to the horizon length $\widetilde{N}$.
\end{assumption}

Note that Assumption \ref{Stability under Adaptation:ass:enhanced stabilizing} is indeed very weak since for one we allow for non--monotone developments of the suboptimality degree $\alpha(\cdot)$ if the horizon length is increased which may occur as shown in \citet{DPM2007}. Here, we only make sure that if a certain suboptimality degree $\overline{\alpha} \in (0, 1)$ holds for a horizon length $N$, then the estimate $\alpha(\widetilde{N})$ does not drop below zero if the horizon length $\widetilde{N}$ is increased.

Considering the value of $l(x, \mu_{\widetilde{N}}(x))$, we notice that it may tend to zero if $\widetilde{N}$ is increased, hence we have that $C_l$ is in general unbounded. The special case $l(x, \mu_{\widetilde{N}}(x)) = 0$, however, states that the equilibrium of our problem has been reached and can be neglected in this context since this implies $l(x, \mu_{N}(x)) = 0$ allowing for arbitrary $C_l$.

Given Assumption \ref{Stability under Adaptation:ass:enhanced stabilizing}, we obtain stability and a performance estimate of the closed loop in the context of changing horizon lengths similar to Proposition \ref{Suboptimality:prop:trajectory a posteriori estimate}.

\begin{theorem}[Stability of Adaptive MPC]\label{Stability under Adaptation:thm:stability of adaptive mpc}
	Consider $\overline{\alpha} \in (0, 1)$ and a sequence $(N_i)_{i \in \N_0}$, $N_i \in \N$, where $N^\star = \max \{ N_i \, | \; i \in \N \}$, such that the MPC feedback law $\mu_{(N_i)}$ defining the closed loop solution \eqref{Setup:closed loop solution} guarantees
	\begin{align}
		\label{Stability under Adaptation:thm:stability of adaptive mpc:eq1}
		V_{N_i}(x(i)) \geq V_{N_i}(x(i+1)) + \overline{\alpha} l(x(i), \mu_{N_i}(x(i)))
	\end{align}
	for all $i \in \N_0$. If additionally Assumption \ref{Stability under Adaptation:ass:enhanced stabilizing} is satisfied for all pairs of initial values and horizons $(x(i), N_i)$, $i \in \N_0$, then we obtain
	\begin{align*}
		\alpha_{C} V_\infty(x(n)) \leq \alpha_{C} V_\infty^{\mu_{(N_i)}}(x(n)) \leq V_{N^\star}(x(n)) \leq V_\infty(x(n))
	\end{align*}
	to hold for all $n \in \N_0$ where $\alpha_{C} := \min\limits_{i \in \N_{\geq n}} C_l^{(i)} C_\alpha^{(i)} \overline{\alpha}$.
\end{theorem}
\textbf{Proof:}
	Given a pair $(x(i), N_i)$, Assumption \ref{Stability under Adaptation:ass:enhanced stabilizing} guarantees $\alpha(N_i) \leq \alpha(\widetilde{N}) / C_\alpha^{(i)}$ for $\widetilde{N} \geq N_i$. Now we choose $\widetilde{N} = N^\star$ within this local suboptimality estimation. Hence, we obtain $\overline{\alpha} \leq \alpha(N_i) \leq \alpha(N^\star) / C_\alpha^{(i)}$ using the relaxed Lyapunov inequality \eqref{Stability under Adaptation:thm:stability of adaptive mpc:eq1}. Multiplying by the stage cost $l(x(i), \mu_{N_i}(x(i)))$, we can conclude
	\begin{align*}
		& \overline{\alpha} l(x(i), \mu_{N_i}(x(i))) \leq  \\
		& \leq \frac{\alpha(N^\star)}{C_\alpha^{(i)}} l(x(i), \mu_{N_i}(x(i))) \\
		& = \frac{V_{N^\star}(x(i)) - V_{N^\star}(f(x, \mu_{N^\star}(x)))}{C_\alpha^{(i)} l(x(i), \mu_{N^\star}(x(i)))} l(x(i), \mu_{N_i}(x(i))) \\
		& \leq \frac{V_{N^\star}(x(i)) - V_{N^\star}(f(x, \mu_{N_i}(x)))}{C_\alpha^{(i)} C_l^{(i)}}
	\end{align*}
	using \eqref{Stability under Adaptation:thm:stability of adaptive mpc:eq1} and Assumption \ref{Stability under Adaptation:ass:enhanced stabilizing}. Summing the running costs along the closed loop trajectory reveals
	\begin{align*}
		\alpha_{C} \sum\limits_{i = n}^{K} l(x(i), \mu_{N_i}(x(i))) \leq V_{N^\star}(x(n)) - V_{N^\star}(x(K+1))
	\end{align*}
	where we defined $\alpha_{C} := \min\limits_{i \in \N_{\geq n}} C_l^{(i)} C_\alpha^{(i)} \overline{\alpha}$. \\
	Since $V_{N^\star}(x(K+1)) \geq 0$ holds, taking $K$ to infinity reveals
	\begin{align*}
		\alpha_{C} V_\infty^{\mu_{(N_i)}}(x(n)) & = \alpha_C \lim\limits_{K \rightarrow \infty} \sum\limits_{i = n}^{K} l(x(i), \mu_{N_i}(x(i))) \\
		& \leq V_{N^\star}(x(n))
	\end{align*}
	Since the $\alpha V_\infty (x(n)) \leq \alpha_{C} V_\infty^{\mu_{(N_i)}}(x(n))$ and $V_{N^\star}(x(n)) \leq V_\infty(x(n))$ hold by the principle of optimality, the assertion follows. \qed

Comparing Proposition \ref{Suboptimality:prop:trajectory a posteriori estimate} and Theorem \ref{Stability under Adaptation:thm:stability of adaptive mpc}, we see that the closed loop estimate $\alpha_C$ may be smaller than the local suboptimality bound $\overline{\alpha}$ but due to $C_l, C_\alpha > 0$ we can guarantee $\alpha_C > 0$. Yet, $\alpha_C$ may become very small depending on $C_\alpha$ and $C_l$ from Assumption \ref{Stability under Adaptation:ass:enhanced stabilizing}.

%%%%%%%%%%%%%%%%%%%%%%%%%%%%%%%%%%%%%%%%%%%%%%%%%%%%%%%%%%%%%%%%%%%%%%%%%%%%%%%%
\section{A simple Adaptation Strategy}
\label{Section:A simple Adaptation Strategy}
%%%%%%%%%%%%%%%%%%%%%%%%%%%%%%%%%%%%%%%%%%%%%%%%%%%%%%%%%%%%%%%%%%%%%%%%%%%%%%%%

Since now we have shown asymptotic stability of a MPC closed loop trajectory with varying optimization horizon, we show a very simple approach to guarantee the local suboptimality requirement $\alpha(N_i) \geq \overline{\alpha}$. To this end, we assume the system to be controllable, i.e.
\begin{assumption}\label{A simple Adaptation Strategy:ass:stabilizable}
	Given $\overline{\alpha} \in (0, 1)$, for all $x_0 \in \X$ there exists a finite horizon length $\overline{N} = N(x_0) \in \N$ such that the relaxed Lyapunov inequality \eqref{Suboptimality:prop:trajectory a posteriori estimate:eq1} holds with $\alpha(N) \geq \overline{\alpha}$.
\end{assumption}
\begin{theorem}[Shortening Strategy]\label{A simple Adaptation Strategy:thm:stepsize shortening}
	Consider an optimal control problem \eqref{Setup:open loop control}, \eqref{Setup:open loop solution} with initial value $x(i)$, horizon $N_i \in \N$ and fixed suboptimality bound $\overline{\alpha} \in (0, 1)$ and denote the optimal control sequence by $u^\star$. Suppose there exists an integer $\overline{i} \in \N_0$, $0 \leq \overline{i} < N_i$ such that
	\begin{align}
		& V_{N_i - k}(x_{u_N}(k, x(i))) - V_{N_i - k}(x_{u^\star}(k + 1, x(i))) \geq \nonumber \\
		\label{A simple Adaptation Strategy:thm:stepsize shortening:eq1}
		& \geq \overline{\alpha} l(x_{u^\star}(k, x(i)), \mu_{N_i - k}(x_{u^\star}(k, x(i))))
	\end{align}
	holds true for all $0 \leq k \leq \overline{i}$. Then, setting $N_{i + k}= N_i - k$ and $\mu_{N_{i + k}}(x(i + k))= u^\star(k)$ for $0 \leq k \leq \overline{i} - 1$, inequality \eqref{Stability under Adaptation:thm:stability of adaptive mpc:eq1} holds for $k = i, \ldots, i + \overline{i} - 1$ with $\alpha = \overline{\alpha}$.
\end{theorem}
\textbf{Proof:}
	The proof follows directly from the fact that for $\mu_{N_{i + k}}(x(i + k))= u^\star(k)$ the closed loop trajectory satisfies $x(i + k) = x_{u^\star}(k, x(i))$. Hence, \eqref{Stability under Adaptation:thm:stability of adaptive mpc:eq1} follows from \eqref{A simple Adaptation Strategy:thm:stepsize shortening:eq1}. 

With the choice $N_{i+k}=N_i-k$, due to the principle of optimality we obtain that the optimal control problems within the next $\overline{i} - 1$ NMPC iterations are already solved since $\mu_{N_{i+k}}(x(i+k))$ can be obtained from the optimal control sequence $u^\star(\cdot) \in \U^N(x(i))$ computed at time $i$. This implies that the most efficient way for the reducing strategy is not to reduce $N_i$ itself but rather to reduce the horizons $N_{i + k}$ by $k$ for the subsequent sampling instants $i + 1, \ldots, i + \overline i$, i.e., we choose the initial guess of the horizon $N_{i+1} = N_i-1$. Still, if the a posteriori estimate is used, the evaluation of \eqref{A simple Adaptation Strategy:thm:stepsize shortening:eq1} requires the solution of an additional optimal control problem in each step. \\
In order to to use the \textit{a priori} estimate given by Theorem \ref{Suboptimality:thm:apriori} the following result can be used as a shortening strategy:

\begin{theorem}[A priori Shortening Strategy]\label{A simple Adaptation Strategy:thm:stepsize shortening2}
	Consider a optimal control problem \eqref{Setup:open loop control}, \eqref{Setup:open loop solution} with initial value $x(i)$ and horizon $N_i, \hat{N} \in \N$, $N_i \geq \hat{N} \geq 2$ and denote the optimal control sequence by $u^\star$. Moreover, the suboptimality bound $\overline{\alpha} \in (0, 1)$ is fixed inducing some $\overline{\gamma}(\cdot)$ via \eqref{Suboptimality:thm:apriori:eq1}. Suppose there exists an integer $\overline{i} \in \N_0$, $0 \leq \overline{i} < N_i - N_0 - 1$ such that for all $0 \leq k \leq \overline{i}$ there exist $\gamma_{i} < \overline{\gamma}(N_i - k)$ satisfying
	\begin{align}
		\label{A simple Adaptation Strategy:thm:stepsize shortening2:eq1}
		& \frac{V_{N_0}(x_{u^\star}(N_i - \hat{N}, x(i)))}{\gamma_{i} + 1} \leq \\
		& \leq \max_{j = 2, \ldots, \hat{N}} l(x_{u^\star}(N_i - j, x(i)), \mu_{j - 1}(x_{u^\star}(N_i - j, x(i)))) \nonumber \\
		\label{A simple Adaptation Strategy:thm:stepsize shortening2:eq2}
		& \frac{V_{k_i}(x_{u^\star}(N_i - k_{i}, x(i)))}{\gamma_{i} + 1} \leq \\
		& \leq l(x_{u^\star}(N_i - k_{i}, x(i)), \mu_{k_{i}}(x_{u^\star}(N - k_{i}, x(i)))) \nonumber
	\end{align}
	for all $k_{i} \in \{\hat{N} + 1, \ldots, N_i - k\}$. Then, setting $N_{i + k}= N_i - k$ and $\mu_{N_{i + k}}(x(i + k))= u^\star(k)$ for $0 \leq k \leq \overline{i} - 1$, inequality \eqref{Stability under Adaptation:thm:stability of adaptive mpc:eq1} holds for $k = i, \ldots, i + \overline{i} - 1$ with $\alpha = \overline{\alpha}$.
\end{theorem}
\textbf{Proof:}
	Since \eqref{A simple Adaptation Strategy:thm:stepsize shortening2:eq1}, \eqref{A simple Adaptation Strategy:thm:stepsize shortening2:eq2} hold for $k = 0$, Theorem \ref{Suboptimality:thm:apriori} guarantees that the local suboptimality degree is at least as large as $\overline{\alpha}$. If $\overline{i} > 0$ holds, we can make use of the fact that for $\mu_{N_{i + k}}(x(i + k))= u^\star(k)$ the closed loop trajectory satisfies $x(i + k) = x_{u^\star}(k, x(i))$. By \eqref{A simple Adaptation Strategy:thm:stepsize shortening2:eq1}, \eqref{A simple Adaptation Strategy:thm:stepsize shortening2:eq2}, we obtain Assumption \ref{Stability under Adaptation:ass:enhanced stabilizing} to hold for $k = i, \ldots, i + \overline{i} - 1$. Accordingly, the assertion follows from Theorem \ref{Suboptimality:thm:apriori} which concludes the proof.\qed

Note that while the a priori estimate from Theorem \ref{Suboptimality:thm:apriori} is slightly more conservative than the result from Proposition \ref{Suboptimality:prop:trajectory a posteriori estimate}, it is also computationally less demanding if the value $N_0$ is small.\\
In contrast to this efficient and simple shortening strategy, it is quite difficult to obtain efficient methods for prolongating the optimization horizon $N_i$.
In order to obtain a simple prolongating strategy, we invert the approach of Theorem \ref{A simple Adaptation Strategy:thm:stepsize shortening}, i.e. we iteratively increase the parameter $N$ until the requirement $\alpha(N_i) \geq \overline{\alpha}$ is satisfied.

\begin{theorem}[Prolongation Strategy]\label{A simple Adaptation Strategy:thm:stepsize prolongation}
	Consider an optimal control problem \eqref{Setup:open loop control}, \eqref{Setup:open loop solution} with initial value $x(i)$ and $N_i \in \N$. Moreover, for fixed $\overline{\alpha} \in (0, 1)$ suppose Assumption \ref{A simple Adaptation Strategy:ass:stabilizable} to hold. Then, any algorithm which iteratively increases the optimization horizon $N_i$ terminates in finite time and computes a horizon length $N_i$ such that \eqref{Stability under Adaptation:thm:stability of adaptive mpc:eq1} holds with local suboptimality degree $\overline{\alpha}$.
\end{theorem}
\textbf{Proof:}
	Follows directly from Assumption \ref{A simple Adaptation Strategy:ass:stabilizable}.\qed

Unfortunately, if \eqref{Stability under Adaptation:thm:stability of adaptive mpc:eq1} does not hold, it is in general difficult to assess by how much $N_i$ should be increased such that \eqref{Stability under Adaptation:thm:stability of adaptive mpc:eq1} holds for the increased $N_i$. The most simple strategy of increasing $N_i$ by one in each iteration shows satisfactory results in practice, however, when starting the iteration with $N_i$, in the worst case it requires us to check \eqref{Stability under Adaptation:thm:stability of adaptive mpc:eq1} $\overline{N} - N_n + 1$ times at each sampling instant. In contrast to the shortening strategy, the principle of optimality cannot be used here to establish a relation between the optimal control problems for different $N_i$ and, moreover, these problems may exhibit different solution structures which makes it a hard task to provide a suitable initial guess for the optimization algorithm. 

%%%%%%%%%%%%%%%%%%%%%%%%%%%%%%%%%%%%%%%%%%%%%%%%%%%%%%%%%%%%%%%%%%%%%%%%%%%%%%%%
\section{Numerical Results}
\label{Section:Numerical Results}
%%%%%%%%%%%%%%%%%%%%%%%%%%%%%%%%%%%%%%%%%%%%%%%%%%%%%%%%%%%%%%%%%%%%%%%%%%%%%%%%

To illustrate the effects of using an adaptive NMPC, we consider a highrack warehouse
\begin{align}
	\ddot{\chi}(t) & = u_1(t) \nonumber \\
	\ddot{\upsilon}(t) & = u_2(t) \\
	\ddot{\phi}(t) & = -k\dot{\phi}(t) - \frac{g}{\upsilon}(t) \sin(\phi(t)) - u_1(t) \cos(\phi(t)) \nonumber
\end{align}
where for simplicity of exposition the rope is modeled as a pendulum with variable length. Here, $\chi$ denotes the position of the crab along the highrack, $\upsilon$ represents the length of the rope of the crane and $\phi$ corresponds to the angle of deflection of the rope. Moreover, $g = 9.81$ and $k = 0.1$ denote the gravitational constant and the inertia of the angle of the rope, respectively.

For this example, we use MPC to generate a feedback for a representative transport action of a pallet from $\chi_0 = -3$, $\upsilon_0 = 5$ to $\chi_{\text{ref}} = 3$, $\upsilon_{\text{ref}} = 2$ (with zero derivatives in initial and target position) while maintaining the state and control constraints $\X = [-5, 5]^2 \times [1, 4] \times [-1, 2] \times [-1, 1] \times \R$ and $\U = [-5, 5] \times [-1, 2]$. To this end, we use the running cost
\begin{align*}
	l(x,u) = \int\limits_0^T & c_1 \dot{\phi}^2(t) \upsilon^2(t)  + c_2 g \upsilon(t) ( 1 - \cos(\phi(t)) ) \\[-4mm]
		& + c_3 (\chi(t) - \hat{\chi})^2 + c_4 \dot{\chi}^2(t) + c_5 (\upsilon(t) - \hat{\upsilon})^2 \\
		& + c_6 \dot{\upsilon}^2(t) + c_7 \left( u_1^2(t) + u_2^2(t) \right) dt
\end{align*}
with constants $c_1 = 0.25$, $c_2 = 0.5$, $c_3 = 40$, $c_4 = c_5 = c_6 = 20$ and $c_7 = 0.1$ and the sampling period $T = 0.2$. To solve the optimal control problem arising throughout the MPC procedure, we use a direct approach, i.e. discretize the continuous time problem and use an SQP method to solve the resulting optimization problem. Here, we set the tolerance levels $\text{tol}_{\text{ODE}} = 10^{-9}$ and $\text{tol}_{\text{SQP}} = 10^{-6}$ for the differential equation solver and the optimization method respectively.

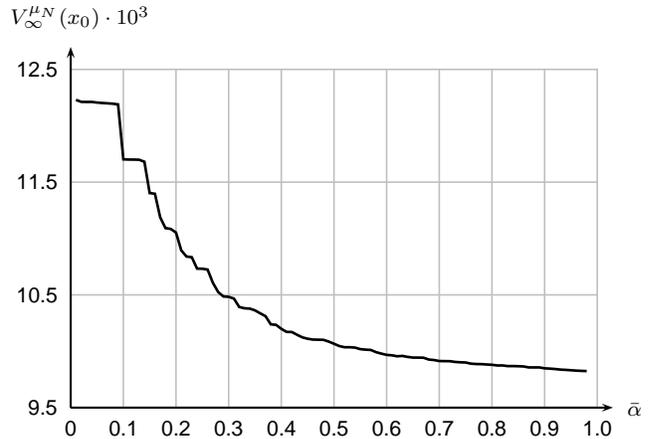
\begin{figure}[ht!]
\begin{center}
\begin{pdfpic}
	\psset{xunit=7cm,yunit=1.5cm,runit=1cm}
	\sffamily
	\def\pshlabel#1{\sffamily\footnotesize #1}
	\def\psvlabel#1{\sffamily\footnotesize #1}
	\begin{pspicture}(-0.12,-0.2)(1.1,3.6)
		\psgrid[yunit=1.0,xunit=0.1,subgriddiv=1,gridcolor=lightgray,gridlabels=0,gridwidth=0.6pt](0,0)(10,3)
		\psaxes[Oy=9.5,Dx=0.1,Dy=1.0,tickstyle=bottom,ticksize=0mm]{->}(1.03,3.2)
		\put(7.4,-0.1){\footnotesize{$\bar{\alpha}$}}
		\put(-0.8,5.1){\footnotesize{$V_\infty^{\mu_N} (x_0) \cdot 10^{3}$}}
		\psline[linecolor=black,linewidth=1.0pt](0.01,2.72984)(0.02,2.71224)(0.03,2.71145)(0.04,2.71175)(0.05,2.70591)(0.06,2.7021)(0.07,2.69917)(0.08,2.69577)(0.09,2.68954)(0.1,2.20244)(0.11,2.2004)(0.12,2.20007)(0.13,2.19872)(0.14,2.18109)(0.15,1.90315)(0.16,1.89709)(0.17,1.68794)(0.18,1.5925)(0.19,1.58508)(0.2,1.55296)(0.21,1.39482)(0.22,1.34011)(0.23,1.33449)(0.24,1.23292)(0.25,1.23185)(0.26,1.22469)(0.27,1.10744)(0.28,1.02583)(0.29,0.985983)(0.3,0.982898)(0.31,0.967358)(0.32,0.892798)(0.33,0.881194)(0.34,0.878421)(0.35,0.86274)(0.36,0.836169)(0.37,0.811576)(0.38,0.738571)(0.39,0.734305)(0.4,0.699037)(0.41,0.671772)(0.42,0.670592)(0.43,0.645912)(0.44,0.623263)(0.45,0.611095)(0.46,0.603729)(0.47,0.602497)(0.48,0.601422)(0.49,0.585912)(0.5,0.566657)(0.51,0.546491)(0.52,0.535822)(0.53,0.535731)(0.54,0.532438)(0.55,0.518943)(0.56,0.514327)(0.57,0.511628)(0.58,0.491927)(0.59,0.479089)(0.6,0.467972)(0.61,0.464103)(0.62,0.455598)(0.63,0.458576)(0.64,0.449939)(0.65,0.443905)(0.66,0.44385)(0.67,0.442741)(0.68,0.426646)(0.69,0.422172)(0.7,0.413295)(0.71,0.412483)(0.72,0.411911)(0.73,0.405432)(0.74,0.40216)(0.75,0.400388)(0.76,0.390231)(0.77,0.386737)(0.78,0.385914)(0.79,0.382618)(0.8,0.380156)(0.81,0.373909)(0.82,0.374274)(0.83,0.368566)(0.84,0.36859)(0.85,0.367402)(0.86,0.365297)(0.87,0.356816)(0.88,0.356626)(0.89,0.356724)(0.9,0.348868)(0.91,0.346065)(0.92,0.342601)(0.93,0.337491)(0.94,0.334998)(0.95,0.331709)(0.96,0.328346)(0.97,0.326103)(0.98,0.324169)
	\end{pspicture}
\end{pdfpic}
\caption{Development of $V_\infty^{\mu_N} (x_0)$ for different suboptimality bounds $\overline{\alpha}$}
\label{Numerical Results:Figure:Vinfty}
\end{center}
\label{fig:V}
\end{figure}

Since the adaptive MPC algorithm allows us to set the lower bound of the degree of suboptimality $\overline{\alpha}$ directly, we first investigate the $\bar{\alpha}$--depending quality of a controlsequence on the closed loop cost $V_\infty^{\mu_N}(x_0)$. To this end, we terminate the algorithm when the condition $l(x(t),u(t)) < 10^{-3}$ is satisfied. The data we obtained for this setting is displayed in Figure \ref{Numerical Results:Figure:Vinfty}.

Here, one can nicely observe that the closed loop costs caused by the adaptive MPC feedback $\mu_{(N_i)}$ are decreasing as the lower bound $\overline{\alpha}$ is enlarged. This is right the behavior one would expect from the theoretical construction. However, using the adaptive MPC approach, larger $\overline{\alpha}$--values not only provides a much better control sequence in terms of generated costs. We also like to mention that the total simulation time required to satisfy the termination criterion is also decreasing as $\overline{\alpha}$ is enlarged which is due to the use of larger optimization horizons $N_i$ throughout the run of the simulation. 

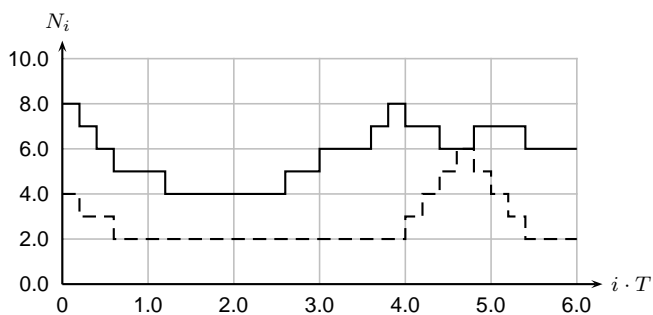
\begin{figure}[ht!]
\begin{center}
\begin{pdfpic}\small
	\psset{xunit=1.14cm,yunit=0.3cm}
	\sffamily
	\def\pshlabel#1{\sffamily\footnotesize #1}
	\def\psvlabel#1{\sffamily\footnotesize #1}
	\begin{pspicture}(-0.1,-0.7)(6,12.5)
		\psgrid[yunit=2.0,xunit=1.0,subgriddiv=1,gridcolor=lightgray,gridlabels=0,gridwidth=0.6pt](0,0)(6,5)
		\psaxes[Oy=0.0,Dx=1.0,Dy=2.0,tickstyle=bottom,ticksize=0mm]{->}(6.3,10.8)
 		\put(7.3,-0.1){\footnotesize{$i \cdot T$}}
 		\put(-0.23,3.4){\footnotesize{$N_i$}}
		\psline[linecolor=black,linestyle=dashed](0,4)(0.2,4)(0.2,3)(0.4,3)(0.4,3)(0.6,3)(0.6,2)(0.8,2)(0.8,2)(1,2)(1,2)(1.2,2)(1.2,2)(1.4,2)(1.4,2)(1.6,2)(1.6,2)(1.8,2)(1.8,2)(2,2)(2,2)(2.2,2)(2.2,2)(2.4,2)(2.4,2)(2.6,2)(2.6,2)(2.8,2)(2.8,2)(3,2)(3,2)(3.2,2)(3.2,2)(3.4,2)(3.4,2)(3.6,2)(3.6,2)(3.8,2)(3.8,2)(4,2)(4,3)(4.2,3)(4.2,4)(4.4,4)(4.4,5)(4.6,5)(4.6,6)(4.8,6)(4.8,5)(5,5)(5,4)(5.2,4)(5.2,3)(5.4,3)(5.4,2)(5.6,2)(5.6,2)(5.8,2)(5.8,2)(6,2)
		\psline[linecolor=black,linestyle=solid](0,8)(0.2,8)(0.2,7)(0.4,7)(0.4,6)(0.6,6)(0.6,5)(0.8,5)(0.8,5)(1,5)(1,5)(1.2,5)(1.2,4)(1.4,4)(1.4,4)(1.6,4)(1.6,4)(1.8,4)(1.8,4)(2,4)(2,4)(2.2,4)(2.2,4)(2.4,4)(2.4,4)(2.6,4)(2.6,5)(2.8,5)(2.8,5)(3,5)(3,6)(3.2,6)(3.2,6)(3.4,6)(3.4,6)(3.6,6)(3.6,7)(3.8,7)(3.8,8)(4,8)(4,7)(4.2,7)(4.2,7)(4.4,7)(4.4,6)(4.6,6)(4.6,6)(4.8,6)(4.8,7)(5,7)(5,7)(5.2,7)(5.2,7)(5.4,7)(5.4,6)(5.6,6)(5.6,6)(5.8,6)(5.8,6)(6,6)
	\end{pspicture}
\end{pdfpic}
\caption{Horizon length $N$ during simulation with $\overline{\alpha}=0.6$(solid) and $\overline{\alpha}=0.2$(dashed)}
\label{Numerical Results:Figure:Ni}
\end{center}
\label{fig:N}
\end{figure}

In Figure \ref{Numerical Results:Figure:Ni}, we additionally plotted the optimization horizon sequences $(N_i)$ for the selected values of $\overline{\alpha}$. This figure demonstrates clearly the horizon incrementations during acceleration-- and deceleration--phases. In particular, a large optimization horizon is required to satisfy the desired decrease in the relaxed Lyapunov inequality \eqref{Stability under Adaptation:thm:stability of adaptive mpc:eq1} upon start of the simulation run which is then reduced as the crab moves towards its destination. In order to reduce the possibly occuring overshoot, the method automatically increases the horizon again. During the final leveling phase, again no large horizons are needed to satisfy \eqref{Stability under Adaptation:thm:stability of adaptive mpc:eq1}.

In Figure \ref{Numerical Results:Figure:Ni}, one can also see that the spike in the horizon length occurs earlier for $\overline{\alpha} = 0.6$. Again, this corresponds to the MPC procedure recognizing the possible overshoot by means of \eqref{Stability under Adaptation:thm:stability of adaptive mpc:eq1}.
\section{Conclusion}
\label{Section:Conclusion}
%%%%%%%%%%%%%%%%%%%%%%%%%%%%%%%%%%%%%%%%%%%%%%%%%%%%%%%%%%%%%%%%%%%%%%%%%%%%%%%%

In this work we have shown stability and suboptimality estimates for model predictive controllers with varying optimization horizon. This result allows for developing strategies to adapt the horizon length instead of using a worst case estimate and control the quality of the feedback directly. \\
Future work may concern reducing the computational effort required to evaluate the suboptimality estimates. Moreover, development and investigation of alternatives to prolongate the optimization horizon will be an issue, i.e. by combining information of several iterates.

% %%%%%%%%%%%%%%%%%%%%%%%%%%%%%%%%%%%%%%%%%%%%%%%%%%%%%%%%%%%%%%%%%%%%%%%%%%%%%%%%
% \bibliography{../bibdatei}
% %%%%%%%%%%%%%%%%%%%%%%%%%%%%%%%%%%%%%%%%%%%%%%%%%%%%%%%%%%%%%%%%%%%%%%%%%%%%%%%%

\end{document}